\documentclass[final]{article}

\usepackage{lmodern}
\usepackage{graphicx}
\usepackage{booktabs}
\usepackage{tikz}
\usepackage{pgfplots}
\pgfplotsset{compat=1.14}
\usepackage{anyfontsize}
\usepackage{hyperref}
\usepackage{caption}
\usepackage{subcaption}
\usepackage{bbm}
\usepackage{amsmath,amsthm,amssymb,amscd}

\renewcommand{\vec}[1]{\mathbf{#1}}

\title{An adaptive Levin method for complicated domains}

\author{Shukui Chen$\mbox{}^{\dagger}$, Kirill Serkh$\mbox{}^{\dagger\,\diamond}$, James Bremer$\mbox{}^{\diamond}$}

\begin{document}

\maketitle

\noindent$\mbox{}^{\dagger}$Department of Computer Science and $\mbox{}^{\diamond}$Department of Mathematics, University of Toronto

\begin{abstract}
\normalsize
\noindent
In this paper we describe an adaptive Levin method for numerically
evaluating integrals of the form $\int_\Omega f(\mathbf x) \exp(i g(\mathbf
x)) \,d\Omega$ over general domains that have been meshed by transfinite
elements. On each element, we apply the multivariate
Levin method over adaptively refined sub-elements, until the integral has
been computed to the desired accuracy. Resonance points on the boundaries of
the elements are handled by the application of the univariate adaptive Levin
method. When the domain does not contain stationary points, the cost of the
resulting method is essentially independent of the frequency, even in the
presence of resonance points.

\end{abstract}

\thispagestyle{empty}

\begin{section}{Introduction}

We consider multivariate oscillatory integrals of the form
\begin{equation}
\label{eqn:oscint}
\int_{\Omega}f(\vec{x})\exp(ig(\vec{x}))\,d\Omega
\end{equation}
where $f$ and $g$ are smooth, slowly varying functions, and $g$ is real-valued. The integration domain $\Omega$ is a bounded region in $\mathbbm{R}^{2}$ with a piecewise smooth boundary. We do not impose any additional conditions on $\nabla g$, and allow for stationary points and resonance points in $\overline{\Omega}$. There exist several approaches for the numerical evaluation of integrals of the form (\ref{eqn:oscint}), including Filon-type methods \cite{gao_extended_2023}, Levin-type methods \cite{olver_quadrature_2006} and numerical steepest descent \cite{huybrechs_construction_2007}.

We demonstrate an adaptive Levin method for numerically evaluating the integral (\ref{eqn:oscint}) over general domains that have been meshed by transfinite triangular elements. By computing the integral as the sum of contributions from the triangular elements, our method can handle arbitrarily complicated domains, and can be run in parallel over the elements. Over each triangular element, we apply the multivariate Levin method on adaptively refined sub-elements of the triangular element, until the integral has been computed to a desired accuracy. Resonance points on the boundaries of the elements are handled by the application of the univariate adaptive Levin method.

\end{section}

\begin{section}{The Levin PDE}

The multivariate Levin method operates by solving the partial differential equation
\begin{equation}
\label{eqn:levinpde}
\mathcal{L}[\vec{p}]:=\nabla\cdot\vec{p}+i\nabla g\cdot\vec{p}=f
\end{equation}
in order to find a vector field $\vec{p}$ such that
\begin{equation}
\nabla\cdot(\vec{p}(\vec{x})\exp(ig(\vec{x})))=f(\vec{x})\exp(ig(\vec{x})).
\end{equation}

By applying the divergence theorem, the integral (\ref{eqn:oscint}) is transformed into
\begin{equation}
\label{eqn:bdryint}
\int_{\Omega}f(\vec{x})\exp(ig(\vec{x}))\,d\Omega=\int_{\Omega}\mathcal{L}[\vec{p}](\vec{x})\exp(ig(\vec{x}))\,d\Omega=\int_{\Gamma}\vec{n}(\vec{x})\cdot\vec{p}(\vec{x})\exp(ig(\vec{x}))\,d\Gamma,
\end{equation}
where $\Gamma$ is the boundary of the integration domain. The univariate integral over the boundary $\Gamma$ can then be evaluated using the univariate adaptive Levin method.

It was long believed that this approach suffers from ``low-frequency breakdown,'' meaning that the accuracy of the calculated value of the integral deteriorates when the integrand is slowly varying, as in the case of a stationary point of $g$. However, it was proven in the univariate case that if a Chebyshev spectral method is used to discretize the differential equation and the resulting linear system is solved via a truncated singular value decomposition, then no low-frequency breakdown occurs \cite{chen_adaptive_2024}. The proof guarantees that in the univariate case, high accuracy can be obtained regardless of the magnitude of $g'$ and whether or not it has zeros. The same argument can be readily generalized to the multivariate case to show that stationary points of $g$ need not be a restriction when applying the multivariate Levin method adaptively.

\end{section}

\begin{section}{Numerical Discretization}
To solve the Levin PDE numerically, we seek a vector field $\vec{p}$ of the form $\vec{p}(\vec{x})=\sum_{j=0}^{n}p_{j}\varphi_{j}(\vec{x})$ where $\{\varphi_{0}, \dots, \varphi_{n}\}$ is the set of all monomials of degree less than or equal to $k$, so that $n+1=\text{dim}(\mathcal{P}_{k})$. Using the collocation method, we solve the system
\begin{equation}
\label{eqn:linsys}
\sum_{j=0}^{n}\mathcal{L}[\varphi_{j}](\vec{x}_{k})p_{j}=f(\vec{x}_{k}), \quad k=0, \dots, m,
\end{equation}
for a set of collocation nodes $\{\vec{x}_{0}, \dots, \vec{x}_{m}\}$, where $m+1=\text{dim}(\mathcal{P}_{\ell})$, see \cite{vioreanu_spectra_2014}. The resulting linear system is then solved using a truncated singular value decomposition. Consistent with our analysis in the univariate case, the discretized PDE can be solved with better accuracy if the collocation nodes have slightly higher order than the basis, so that $\ell > k$.

By incorporating phase function methods for ordinary differential equations, the multivariate Levin method can be applied to a large class of oscillatory integrals involving special functions, including products of such functions and compositions of such functions with slowly-varying functions, without the need for any symbolic or arbitrary precision computations \cite{chen_adaptive_2024}.

\end{section}

\begin{section}{Resonance Points}
A point of resonance is any point $\vec{x}\in\Gamma$ such that the vector $\nabla g(\vec{x})$ is orthogonal to the boundary \cite{iserles_computation_2006}. The resonance points represent stationary points of the oscillator in the boundary integral (\ref{eqn:bdryint}), which means that the oscillatory boundary integral cannot be efficiently evaluated using the non-adaptive univariate Levin method of fixed order without sampling at the stationary points. A nonresonance condition is usually imposed to require the vector $\nabla g(\vec{x})$, $\vec{x}\in\Gamma$, to be nowhere orthogonal to the boundary.

We note that the vector field $\vec{p}$ can often be found to high accuracy, in the sense that the error $\|\mathcal{L}[\vec{p}]-f\|_{L^{\infty}(\Omega)}$ is small regardless of whether or not resonance points are present. This suggests that the integral (\ref{eqn:oscint}) can be accurately approximated by the boundary integrals (\ref{eqn:bdryint}). Therefore, we apply the univariate adaptive Levin method to efficiently evaluate the boundary integrals, since it can handle the case of stationary points automatically. As a consequence, we can evaluate of the integral (\ref{eqn:oscint}) efficiently by applying the Levin method adaptively and recursively, regardless of the magnitude of $\nabla g$ and whether stationary points or resonance points are found in $\overline{\Omega}$.
\end{section}

\begin{section}{Numerical experiment}
We demonstrate the efficiency and accuracy of our method by computing the integral
\begin{equation}
\label{eqn:example}
\int_{\Omega}\nabla\cdot(G(\vec{x},\vec{x}')\nabla u(\vec{x}'))\,d\vec{x}', \quad \vec{x}=\vec{0},
\end{equation}
where $G(\vec{x},\vec{x}')=\frac{i}{4}H_{0}^{(1)}(\omega\|\vec{x}-\vec{x}'\|)$ is the Green's function of the 2D Helmholtz equation and $u(\vec{x}')=\|\vec{x}'\|_{2}^{2}$ is chosen to be a smooth function. We use the non-oscillatory phase function for Bessel's equation to transform the integral (\ref{eqn:example}) into the form (\ref{eqn:oscint}).

The integration domain $\Omega$ (Figure \ref{fig:real_part}) was constructed and tessellated using Gmsh and its built-in CAD kernel \cite{geuzaine_gmsh_2009}. After tessellation, we constructed piecewise Chebyshev expansions to represent the boundary as curved segments, and built the geometric mapping using transfinite interpolation \cite{perronnet_interpolation_1998}.

\begin{figure}
\centering
\includegraphics[width=\textwidth]{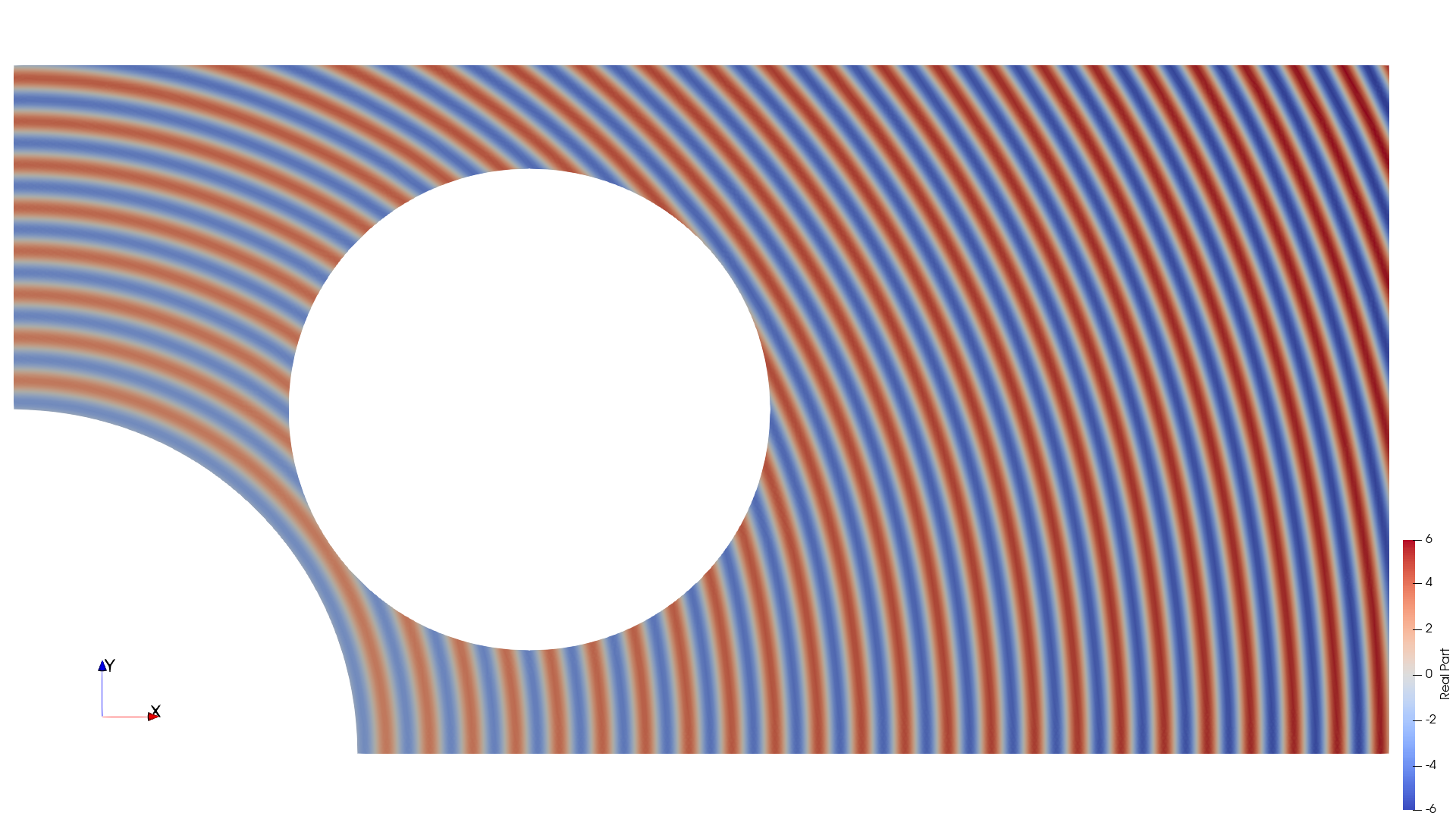}
\caption{Real part of the integrand of (\ref{eqn:example}) with $\omega=100$.
\label{fig:real_part}
}
\end{figure}

We point out that, due to the geometry, every point on the curved boundary in the bottom left is a resonance point. Moreover, due to the tessellation, resonance points are also found on some edges of the interior triangles. However, our method was nonetheless able to evaluate the integral (\ref{eqn:example}) in a nearly frequency-independent manner.

For estimating the absolute error obtained by our method, we transformed the integral into a boundary integral by noting that the integrand is the divergence of a vector field. We then evaluated the boundary integral adaptively using a 30-point Gauss-Legendre quadrature rule.

\end{section}

\begin{section}{Numerical experiment}

We present in Figure~\ref{fig:plots} the results of the numerical
experiments computing the integral (\ref{eqn:example}), in order to
illustrate the properties of the multivariate adaptive Levin method. The
code for the experiment was written in Fortran and compiled with version
13.2.0 of the GNU Fortran compiler. The experiments were performed on a
laptop computer equipped with an Intel i9-12900HK processor and 32GB of
memory. We used OpenMP to parallize over the triangular elements.

\begin{figure}
\centering
\begin{subfigure}[b]{0.65\textwidth}
\centering
\includegraphics[width=\textwidth]{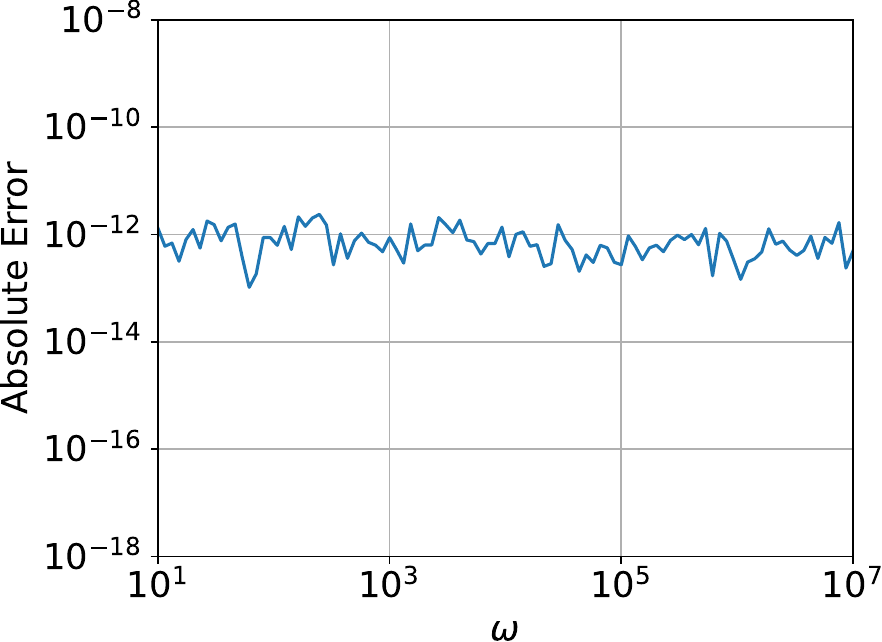}
\end{subfigure}

\begin{subfigure}[b]{0.65\textwidth}
\centering
\includegraphics[width=\textwidth]{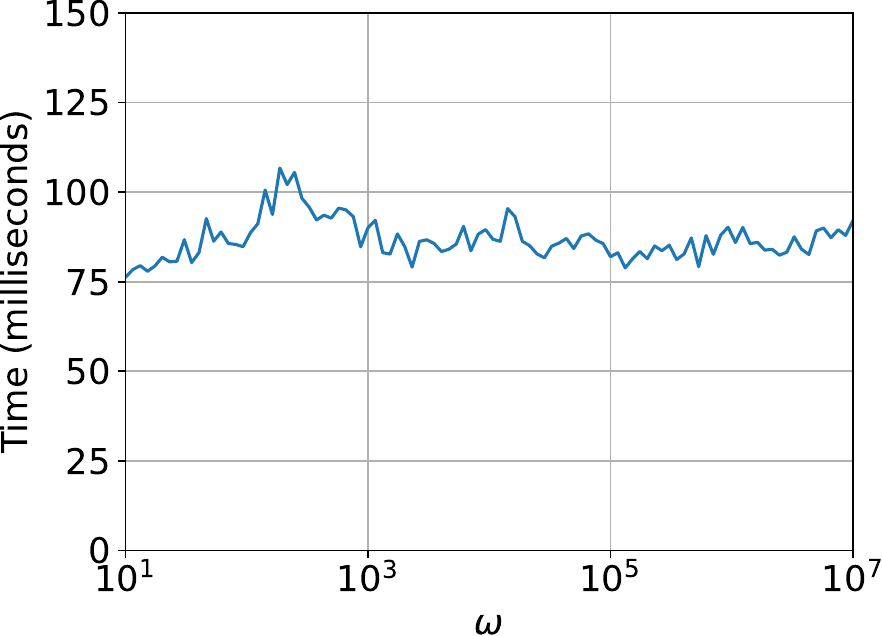}
\end{subfigure}
\caption{The results of the numerical experiments computing the integral (\ref{eqn:example}). The top plot gives the error in the value of the integral computed via the multivariate adaptive Levin method as a function of $\omega$. The bottom plot gives the running time in milliseconds as a function of $\omega$.
\label{fig:plots} }
\end{figure}

\end{section}

\newpage 

\begin{section}{References}
\footnotesize{\bibliographystyle{ieeetr}\bibliography{poster2}}
\end{section}

\end{document}